\documentclass[12pt]{iopart}
\usepackage{epsfig,graphicx,graphics,cite}
\usepackage{iopams}
\begin{document}

\newtheorem{theorem}{Theorem}
\newtheorem{definition}[theorem]{Definition}
\newtheorem{corollary}[theorem]{Corollary}
\newtheorem{remark}[theorem]{Remark}
\title{Exact inversion of the conical Radon transform with a fixed opening angle}
\author{Rim Gouia-Zarrad$^1$ and Gaik Ambartsoumian$^2$}
\address{$^1$ Department of Mathematics and Statistics, American University of Sharjah, PO Box 26666, Sharjah, UAE}
\address{$^2$ Department of Mathematics, University of Texas at Arlington, Arlington, TX 76017-0408, USA}
\ead{\mailto{rgouia@aus.edu} and \mailto{gambarts@uta.edu}}

\begin{abstract}
We study a new class of Radon transforms defined on circular cones called the conical Radon transform. In  $\mathbb{R}^3$  it maps a function to its surface integrals
over circular cones, and in $\mathbb{R}^2$ it maps a function to its integrals along two rays with a common vertex. Such transforms appear in various mathematical models arising in medical imaging, nuclear industry and homeland security. This paper contains new results about inversion of conical Radon transform with fixed opening angle and vertical central axis in $\mathbb{R}^2$ and $\mathbb{R}^3$. New simple explicit inversion formulae are presented in these cases. Numerical simulations were performed to demonstrate the efficiency of the suggested algorithm in 2D.
\end{abstract}
\maketitle

\section{Introduction}

The Conical Radon Transform (CRT) integrates a function over circular cones in $\mathbb{R}^3$, and along coupled rays with a common vertex in $\mathbb{R}^2$. The interest towards such transforms during the last decade was triggered by the connection between the CRT and some mathematical models arising in medical imaging \cite{gaik, gaiknew, basko1, basko2, cree, florescu1, florescu2, florescu3, nguyen5, nguyen1, nguyen3, nguyen6}, homeland security \cite{kuchment1}, nuclear industry and astrophysics \cite{astro}.

The first work concerning CRT was done by Cree and Bones who considered the three dimensional problem in relation to Compton cameras \cite{cree}. They showed that in the case when the vertices of integration cones are restricted to a plane, while their symmetry axis is normal to the plane and the opening angles vary, the CRT can be uniquely inverted.

Basko and collaborators analyzed the related problem for the planar CRT with the vertices of coupled rays restricted to a line and the image function $f$ supported on one side of that line \cite{basko1,basko3}. They showed that the problem of inverting the CRT in this setup can be reduced to the inversion of the Radon transform of another function related to $f$ by a certain symmetry. They extended their idea to the 3D case in \cite{basko2} with the aid of a spherical harmonics expansion, in order to convert the cone-surface integrals into the classical Radon projections.

Some further results for 2D and 3D CRT in these setups were presented by Nguyen, Truong and collaborators. As in the previous works, they considered CRT with vertices restricted to a plane in 3D and a line in 2D, and varying opening angle. We refer the readers to \cite{nguyen5,nguyen1,nguyen3,nguyen6} for further details.

In a series of works \cite{florescu1, florescu2,florescu3} Florescu, Markel, and Schotland investigated image reconstruction problems from CRT appearing in 2D optical tomography. Here the axes of symmetry and the opening angle between the coupled rays of integration were fixed, while the locations of their vertices were not restricted. The authors proved that in this setup CRT can be inverted and provided an exact reconstruction formula.

 Ambartsoumian studied the related problem in a circular setup, where the opening angle of CRT was fixed, one of the coupled rays was normal to the unit circle, while the locations of vertices varied inside the disc. For functions supported sufficiently deep inside the disc he showed \cite{gaik} that the inversion problem of CRT can be reduced to the inversion of the classical Radon transform. In a subsequent work  \cite{gaiknew} Ambartsoumian and Moon showed that CRT in this setup has a unique inversion without any restriction on the support of $f$ and provided an exact inversion formula.

 It must be noted that different authors use different terminology to denote CRT in 2D. Nguyen, Truong et al, as well as Ambartsoumian and Moon call it a V-line transform (since the integration is done along V-shaped trajectories). Florescu et al call this transform a broken-ray transform. The latter term is used by other authors to denote a transform integrating along trajectories made of multiple connected linear intervals, e.g. see \cite{Hub,joonas} and the references there. To avoid this type of confusion in this paper we will follow the first choice and call the 2D CRT a V-line transform.

In this paper we concentrate on the problem of recovering a function from its conical projections in the case of fixed opening angle and a vertical central axis in both 2D and 3D. The geometry and the setup of CRT in $ \mathbb{R}^3$ is different from the ones considered in previous works and we provide new inversion results for this case. In $\mathbb{R}^2$ the CRT that we consider is similar to the one studied by Florescu et al \cite{florescu1, florescu2,florescu3}, but we provide a much simpler inversion formula, which holds for a wider range of opening angles, and demonstrate its efficiency with numerical simulations.

  The rest of this paper is organized as follows. Section 2 is devoted to the introduction of the concept of CRT in $ \mathbb{R}^2$ and $ \mathbb{R}^3$. In section 3, we present new mathematical results on the inversion of the two dimensional CRT and derive an explicit inversion formula. The proof is based on the use of the Fourier transform and the theory of integral equations. In section 4, we study the three dimensional case. New explicit inversion formula is derived for the CRT on a special class of cones. In section 5, numerical simulations are presented to demonstrate the efficiency of the suggested 2D algorithm. Additional remarks are given in the last section with acknowledgments and bibliography.

\section{Formulation of the problem}

\begin{figure}[h]
\begin{center}
{\epsfig{figure=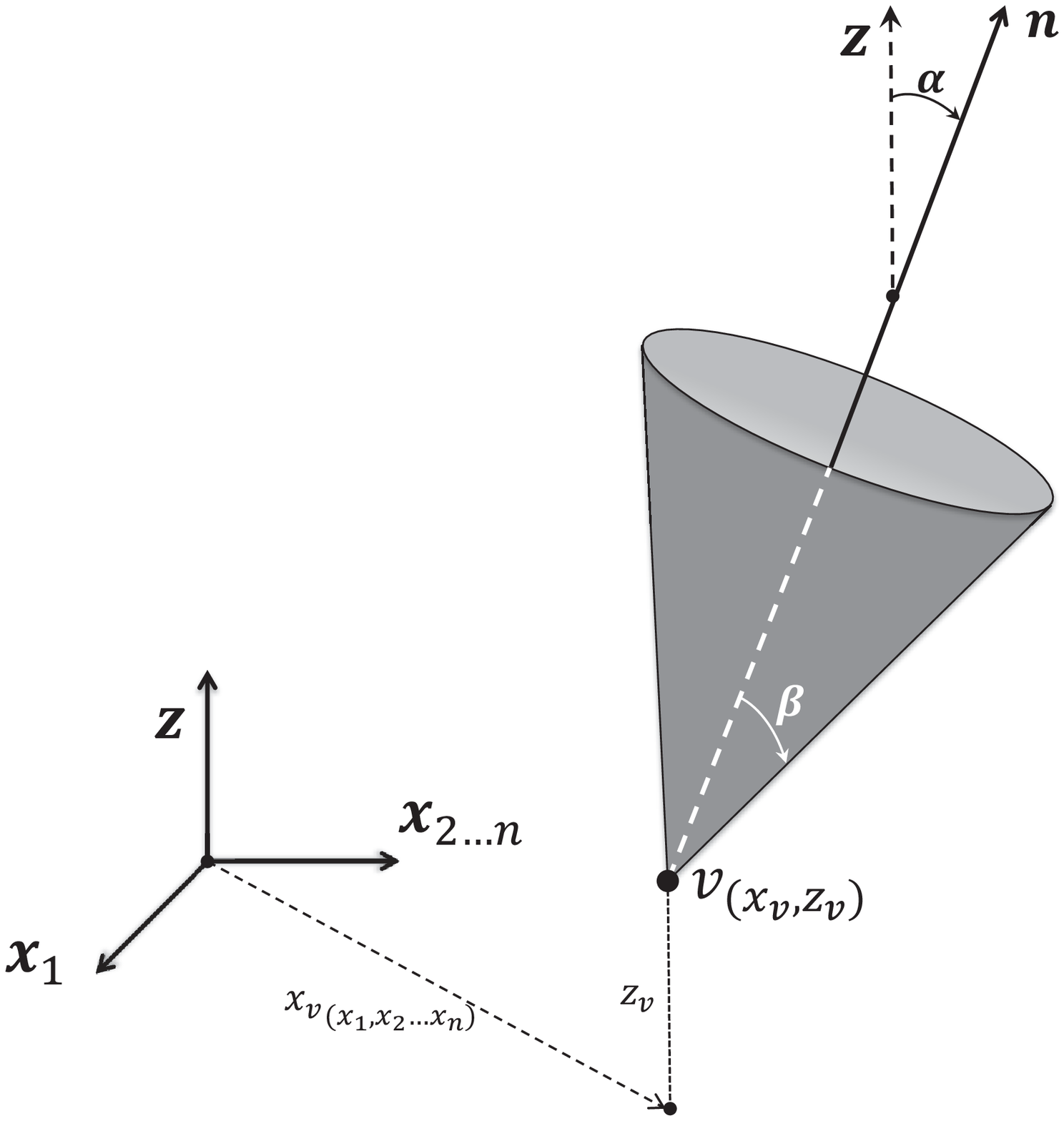,height=3in}}
\end{center}
\vspace{-0.5cm}
\caption{Cone $C(V,\beta, \textbf{n})$.} \label{fig1}
\end{figure}

Consider an $(n+1)$-dimensional $(x,z)$ space, where $n= 2 \mbox{ or } 3$, and $x=(x_1,\ldots,x_n)$. We consider a family of circular cones $C(V,\beta, \textbf{n})$ parameterized by a vertex $V$ of coordinates $(x_v,z_v)$, central axis vector $ \textbf{n}$, and a half opening angle $\beta \in (0,\frac{\pi}{2})$ (see figure \ref{fig1}).\\

If $n=3$ we define the Conical Radon Transform (CRT) of a given function $f(x,z)$ to be the surface integral of the function along a cone $C(V, \beta, \textbf{n})$
\begin{equation}\label{eq1}
g(V,\beta, \textbf{n})=\int_{C(V, \beta, \textbf{n})} f(x,z)\, ds,
\end{equation} \noindent where $ds$ is the surface element on the cone $C(V,\beta, \textbf{n})$.

If $n=2$ the above integration is done along a a pair of rays
with the same vertex $V$ symmetrically located on two sides on the normal vector $ \textbf{n}$. Here $ds$ is an arc length element along the rays, and we call the CRT a V-line Radon transform of $f$.

This paper is dedicated to the problem of recovering the unknown function $f(x,z)$ from the values of an $n$-dimensional family
of its conical projections $g(V,\beta, \textbf{n})$. As with many other generalized Radon transforms with translation or rotation invariance
(e.g. see \cite{rim, gaiknew, florescu1, mypaper,romanov, Kuch, Natt, ref:Norton2D}), our proofs are based on Fourier techniques both for $n=2$ and for $n=3$.\\

\section{Exact inversion formula in 2D}

\begin{figure}[h]
\begin{center}
{\epsfig{figure=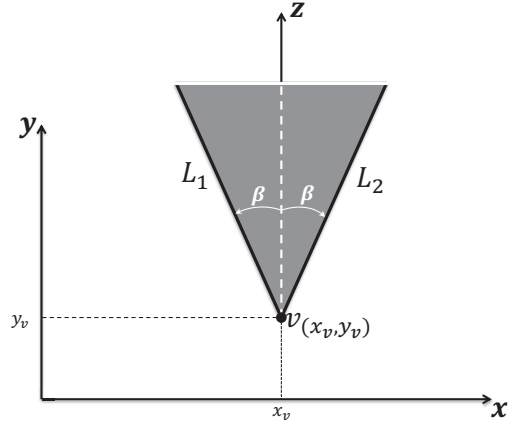,height=3in}}
\end{center}
\vspace{-0.5cm}
\caption{Two-dimensional cone $C((x_v,y_v),\beta, \textbf{y})$.} \label{fig2}
\end{figure}

In two dimensional setting, a circular cone simply consists of two half-lines $L_1$ and $L_2$ with common vertex of coordinates $(x_v,y_v)$ (see figure \ref{fig2}). We will use the existing terminology and call our transform the V-line Radon transform \cite{gaik, gaiknew, nguyen5,nguyen6}. If $g(x_v,y_v,\beta,\textbf{n})$ is known for all possible values of its five arguments then the reconstruction of the function $f(x,y)$ of two variables is an overdetermined problem. To match the dimensions of the data and the function of interest, one can apply additional restrictions on $g$. There are many different ways of reducing the number of parameters, e.g. by considering the vertex $V$ on the $x$-axis, by fixing the half-angle $\beta$ or by fixing the central axis etc. All of these approaches lead to interesting mathematical problems about the invertibility of the V-line Radon transform.


  In this paper, we shall consider the family of cones $C(x_v,y_v,\beta,\textbf{n})$ with fixed half opening angle $\beta$ $ \in (0,\frac{\pi}{2})$, vertical central axis (i.e. parallel to $\textbf{y}$ axis) and no restriction on the vertex. This problem is similar to the one considered in \cite{florescu1}, but we use a different approach to solve it, our inversion formula is much simpler than the one obtained there, and it holds for a wider range of opening angles, since the inversion in \cite{florescu1} holds only for $\beta\in(\pi/4,\pi/2)$.
\begin{theorem}
 Consider a function $f(x,y)\in C^{\infty}_c(\mathcal{D})$, where $ \mathcal{D}=\{(x,y)\in \mathbb{R}^2 \mid 0\leq x \leq x_{max} ,0\leq y \leq y_{max}\}$.
 For $(x_v,y_v)\in \mathbb{R}^2$ we define the V-line Radon transform as follows
\begin{equation}\label{eq2d}
g(x_v,y_v)=\int_{C(x_v,y_v)} f(x,y)\, ds.
\end{equation}
An exact solution of the inversion problem for the V-line Radon transform is given by the formula
\begin{equation}\label{result}f(x,y)=-\frac{\cos{\beta}}{2}\left(\frac{\partial}{\partial y}\, g(x,y)+\tan^2(\beta) \int_{y}^{y_{max}}\frac{\partial^2}{\partial x^2}\,g(x,t)\, dt\right)
\end{equation}

\end{theorem}

\noindent {\bf Proof.}
 The V-line Radon transform is
 \begin{equation}\label{radon}
g(x_v,y_v)=\int_{L_1} f(x,y)\, dl+\int_{L_2} f(x,y)\, dl.
\end{equation}

 \noindent We can represent the equation of the half-lines $L_1$ and $L_2$ in the form
\[
 x=  x_v\pm r \sin(\beta)\,\,\,\,\,\,\,\,\, y= y_v+r\cos(\beta) \]
where $r \geq 0$ is the distance along $L_1$ or $L_2$ measured from the vertex $(x_v,y_v)$.
\noindent Since $\tan(\beta)$ is constant, it is convenient to denote $t=\tan(\beta)$ and $x=x_v\pm (y-y_v)t$. Let us rewrite the equation (\ref{radon}) as \[
g(x_v,y_v)=\int_{y_v}^{y_{max}} \frac{f(x_v+ (y-y_v)t,y)}{\cos(\beta)}\, dy+\int_{y_v}^{y_{max}} \frac{f(x_v- (y-y_v)t,y)}{\cos(\beta)}\, dy.
\]

\noindent The first step consists of applying the Fourier transform with respect to the variable $x_v$ where the Fourier transform $\widehat{g_{\lambda}}( y_v)$ is denoted by

\[
\widehat{g_{\lambda}}( y_v)=\int_{-\infty}^{\infty}g(x_v,y_v)\,e^{-i\lambda x_v}dx_v.
\]
\begin{eqnarray*}
\widehat{g_{\lambda}}( y_v)\cos(\beta)=\int_{-\infty}^{\infty}\int_{y_v}^{y_{max}}f(x_v+ (y-y_v)t,y)\, \,e^{-i\lambda x_v}dy dx_v\\
\hspace{+3cm}+\int_{-\infty}^{\infty}\int_{y_v}^{y_{max}}f(x_v- (y-y_v)t,y)\, \,e^{-i\lambda x_v}dy dx_v\\
\hspace{+3cm}=\int_{y_v}^{y_{max}}\int_{-\infty}^{\infty}f(x_v+ (y-y_v)t,y)\, \,e^{-i\lambda x_v} dx_v dy\\
\hspace{+3cm} +\int_{y_v}^{y_{max}}\int_{-\infty}^{\infty}f(x_v- (y-y_v)t,y)\, \,e^{-i\lambda x_v} dx_vdy
\end{eqnarray*}
\noindent
The equation may be given another form with the following change of variables $X=x_v+ (y-y_v) t$ for the first part and $X=x_v- (y-y_v) t$ for the second part

\begin{eqnarray*}\widehat{g_{\lambda}}( y_v)\cos(\beta)=\int_{y_v}^{y_{max}}\int_{-\infty}^{\infty}f(X,y)\, \,e^{-i\lambda X}e^{i\lambda (y-y_v)t}dXdy \\
\hspace{+3cm} + \int_{y_v}^{y_{max}}\int_{-\infty}^{\infty}f(X,y)\, \,e^{-i\lambda X}e^{-i\lambda (y-y_v)t}dXdy .
\end{eqnarray*}

\noindent This result is further simplified using the Fourier transform \[\widehat{f_{\lambda}}( y)=\int_{-\infty}^{\infty}f(X,y)\,e^{-i\lambda X}dX\] to obtain
\[
\widehat{g_{\lambda}}( y_v)\cos(\beta)=\int_{y_v}^{y_{max}}\widehat{f_{\lambda}}( y)e^{i\lambda (y-y_v)t}dy+\int_{y_v}^{y_{max}}\widehat{f_{\lambda}}( y)e^{-i\lambda (y-y_v)t}dy.
\]

\noindent To simplify the notation we introduce the function $G_{\lambda}( y)$ defined by the following formula:
\[G_{\lambda}( y_v)=\frac{\widehat{g_{\lambda}}( y_v)\cos(\beta)}{2}.\]
\noindent This yields
\begin{equation}\label{G}
G_{\lambda}( y_v)=\int_{y_v}^{y_{max}}\widehat{f_{\lambda}}( y)\cos(\lambda (y-y_v)t)dy.
\end{equation}

\noindent In order to obtain an explicit formula for $\widehat{f_{\lambda}}( y)$, we first differentiate the equation with respect to $y_v$

\begin{equation}\label{G'}
G'_{\lambda}( y_v)=\int_{y_v}^{y_{max}}\widehat{f_{\lambda}}( y)\,\lambda\, t \,\sin(\lambda (y-y_v)t) \,dy -\widehat{f_{\lambda}}( y_v),
\end{equation}
\noindent and combine it with
\begin{equation}\label{G"}
\int_{z}^{y_{max}}G_{\lambda}( y_v)dy_v=\frac{1}{\lambda t}\int_{z}^{y_{max}}\widehat{f_{\lambda}}( y) \,\sin(\lambda (y-z)t) \,dy.
\end{equation}
Finally, we arrive at the following formula
\begin{equation}\label{inversion}
\widehat{ f_{\lambda}}( z)=-G'_{\lambda}( z)+(\lambda \, t)^2\int_{z}^{y_{max}}G_{\lambda}( y_v)dy_v
\end{equation}
By inverse Fourier transform, the solution is
\begin{equation}\label{theo}
f(x,z)=-\frac{\cos{\beta}}{2}\left(\frac{\partial}{\partial z} g(x,z)+\tan^2(\beta) \int_{z}^{y_{max}}\frac{\partial^2}{\partial x^2}g(x,y) dy\right)
\end{equation}

\hfill $\Box$


\vspace{2mm}
\section{Exact inversion formula in 3D}

 \begin{figure}[h]
\begin{center}
{\epsfig{figure=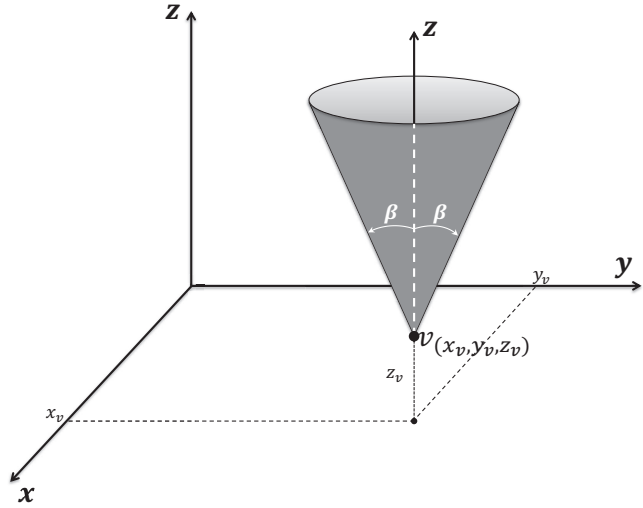,height=3in}}
\end{center}
\vspace{-0.5cm}
\caption{Three-dimensional cone $C((x_v,y_v,z_v),\beta, \textbf{z})$.} \label{fig3}
\end{figure}

In the 3D case of this paper, we consider the family of cones  $C(x_v,y_v,z_v,\beta,\textbf{n})$ with fixed half opening angle $\beta$ $ \in (0,\frac{\pi}{2})$, vertical central axis $\textbf{z}$ and no restriction on the vertex. The problem of integral geometry considered in this section is to reconstruct the unknown function $f(x,y,z)$ in terms of the conical projections $g(x_v,y_v,z_v)$ (see figure \ref{fig3}). To the best of our knowledge, no exact inversion formula is known for this reconstruction problem.

\begin{theorem}
Consider a function $f\in C^{\infty}_c(\mathcal{D})$, where $ \mathcal{D}=\{(x,y,z)\in \mathbb{R}^3 \mid 0\leq x \leq x_{max} ,0\leq y \leq y_{max},0\leq z \leq z_{max}\}$. For $(x_v, y_v, z_v)\in\mathbb{R}^3$ we define the 3D conical Radon transform by
\begin{equation}\label{eq3}
g(x_v,y_v,z_v)=\,\int_{C(x_v,y_v,z_v)}f(x,y,z)\,ds.
\end{equation}
  An exact solution of the inversion problem for CRT is given by \[ \widehat{f}_{\lambda,\mu}(z)=\displaystyle{\frac{\cos{\beta}}{2\pi\tan \beta}\int_{z_{max}}^{z} J_0\left(u(z-x)\right)\left[\frac{d^2}{dx^2}+u^2\right]^2 \int_{z_{max}}^{x}\widehat{g}_{\lambda,\mu}(z_v)\,dz_v\,dx}\]
  where $ \widehat{g}_{\lambda,\mu}(z_v)$ and $\widehat{f}_{\lambda,\mu}(z)$ are the 2D Fourier transforms of the functions $g(x_v,y_v,z_v)$ and
$f(x,y,z)$ with respect to the first two variables.
 \end{theorem}

\noindent {\bf Proof.}
Due to the invariance of the family  $C(x_v,y_v,z_v)$ with respect to translations along the hyperplane $z=0$, the conical Radon transform can be simplified as follows
\[\displaystyle{ g(x_v,y_v,z_v)=\,\int_{C(0,0,z_v)}f(x+x_v,y+y_v,z)ds}.\]
\noindent Then we apply the 2D Fourier transform with respect to the variables $x_v$ and $y_v$ to get
\begin{equation}\label{shift}\displaystyle{ \widehat{g}_{\lambda,\mu}(z_v)=\,\int_{C(0,0,z_v)}\widehat{f}_{\lambda,\mu}(z)e^{i\lambda x}e^{i\mu y}ds}
\end{equation}
\noindent where $ \widehat{g}_{\lambda,\mu}(z_v)$ and $\widehat{f}_{\lambda,\mu}(z)$ are the Fourier
transforms of the functions $g(x_v,y_v,z_v)$ and $f(x,y,z)$.

Let us denote by $S(k,z_v)$ the circles of the intersection of $C(0,0,z_v)$ and the hyperplane $z=k$. We can split the surface element $ds$ into an integration with respect to $dz$ and $dl$, which is the arc length along the circle $S(z,z_v)$.
\[ds=\frac{dl\, dz}{\cos \beta}.\]
\noindent Equation (\ref{shift}) becomes

\begin{equation}\label{V11}
\widehat{g}_{\lambda,\mu}(z_v)=\,\int_{z_v}^{z_{max}}\widehat{f}_{\lambda,\mu}(z)K_{\lambda,\mu}(z,z_v)\,  dz.
\end{equation}

\noindent The kernel is given by the formula
\begin{equation}{\label{kernel}}
\displaystyle K_{\lambda,\mu}(z,z_v)=\,\int_{S_{(z,z_v)}}e^{i(\lambda x+\mu y)}\,\frac{ dl}{\cos{\beta}}.
\end{equation}






\noindent The latter is easily calculated using the equations of the circles $S(k,z_v)$ (we set $x_v=y_v=0$) and introducing the polar coordinates

\[x=r\cos \varphi \,\,\,\,\,\,\,  y=r \sin \varphi \,\,\,\,\,\,\,
 r=(z-z_v)\tan \beta.\]

\noindent We can calculate $dl$ as follows

\[dl=r d\varphi= (z-z_v) \tan \beta\, d\varphi.\]

\noindent Equation (\ref{kernel}) becomes \[\displaystyle K_{\lambda,\mu}(z,z_v)=\,\int_{0}^{2\pi}\frac{e^{i(\lambda x+\mu y)}\, (z-z_v) \tan \beta d\varphi}{\cos{\beta}}\]

\noindent or \[\displaystyle{ K_{\lambda,\mu}(z,z_v)=\frac{(z-z_v) \tan \beta} {\cos{\beta}}\,\int_{0}^{2\pi}e^{i(z-z_v)\tan \beta(\lambda \cos \varphi+\mu \sin \varphi)}\,  d\varphi}.\]

\noindent In the case when $(\lambda,\mu)\ne (0,0)$ let us introduce a new variable $\omega$ as follows
\[\lambda=\sqrt{\lambda^2+\mu ^2}\cos{\omega}, \,\,\,\,\,\,\,\mu=\sqrt{\lambda^2+\mu ^2}\sin{\omega}.\]
Then the expression of the kernel becomes
\[\displaystyle {K_{\lambda,\mu}(z,z_v)=\frac{(z-z_v) \tan \beta} {\cos{\beta}}\,\int_{0}^{2\pi}e^{i(z-z_v)\tan \beta \sqrt{\lambda^2+\mu ^2}(\cos{\omega} \cos \varphi+\sin{\omega} \sin \varphi)}\,  d\varphi},\]
\noindent which yields \[\displaystyle {K_{\lambda,\mu}(z,z_v)=\frac{(z-z_v) \tan \beta} {\cos{\beta}}\,\int_{0}^{2\pi}e^{i(z-z_v)\tan \beta \sqrt{\lambda^2+\mu ^2}(\cos(\omega-\varphi))\,}}  d\varphi.\]

\noindent Using the integral representation of a Bessel function of first kind \[J_0(z)=\frac{1}{2\pi} \,\int_{0}^{2\pi}e^{iz \cos(\theta)} d\theta\] (see \cite{handbook}, page 360), we can express this integral in terms of $J_0$
 \[\displaystyle{K_{\lambda,\mu}(z,z_v)=\frac{2\pi\tan \beta} {\cos{\beta}} (z-z_v) J_0\left((z-z_v)\tan \beta \sqrt{\lambda^2+\mu ^2}\right)}.\] As a result our problem breaks down to the following set of one-dimensional integral equations

\begin{equation}\label{final0}
\displaystyle {\widehat{g}_{\lambda,\mu}(z_v)=\frac{2\pi\tan \beta} {\cos{\beta}}\,\int_{z_v}^{z_{max}}\widehat{f}_{\lambda,\mu}(z)(z-z_v) J_0\left((z-z_v)\tan \beta \sqrt{\lambda^2+\mu ^2}\right)}\,  dz.
\end{equation}
By setting $$\displaystyle {G_{\lambda,\mu}(z_v)=\frac{\cos{\beta}}{2\pi\tan \beta}}\widehat{g}_{\lambda,\mu}(z_v) $$ and $$u=\tan \beta \sqrt{\lambda^2+\mu ^2},$$ one can transform equation (\ref{final0}) to the following expression
\begin{equation}\label{final}
\displaystyle {G_{\lambda,\mu}(z_v)=\,\int_{z_v}^{z_{max}}\widehat{f}_{\lambda,\mu}(z)(z-z_v) J_0\left(u(z-z_v)\right)}\,  dz.
\end{equation}

\noindent We now show how to find $\widehat{f}_{\lambda,\mu}(z)$. For this purpose we introduce the operator ${\mathcal{{H}}}$ defined by
\[\large{\mathcal{{H}}}(F)\equiv\left(\frac{d^2}{dx^2}+u^2\right)F(x)\]
\noindent and apply it to $\displaystyle {\left(\int_{x}^{z_{max}}G_{\lambda,\mu}(z_v)dz_v\right)}$. First, we differentiate our integral using Leibnitz's rule as follows

\begin{eqnarray}\label{diff}
\frac{d^2}{dx^2}\left(\int_{x}^{z_{max}}G_{\lambda,\mu}(z_v)dz_v\right)=\int_{x}^{z_{max}}\widehat{f}_{\lambda,\mu}(z) J_0(u(z-x))  dz \nonumber \\
\hspace{+4.5cm} -u\int_{x}^{z_{max}}\widehat{f}_{\lambda,\mu}(z) (z-x)J_1\left(u(x-z)\right)\,  dz.
\end{eqnarray}


\noindent Here we used the fact that $J'_0(z)=-J_1(z)$ (see \cite{handbook} page 361). We also write

\[\int_{x}^{z_{max}}G_{\lambda,\mu}(z_v)dz_v= \int_{x}^{z_{max}}\,\int_{z_v}^{z_{max}}\widehat{f}_{\lambda,\mu}(z)(z-z_v) J_0(u(z-z_v)) dz\,dz_v, \]
which can be simplified by interchanging the integrals

\begin{equation}\label{bessel1}
\int_{x}^{z_{max}}G_{\lambda,\mu}(z_v)dz_v= \int_{x}^{z_{max}}\widehat{f}_{\lambda,\mu}(z)\,dz \int_{x}^{z}(z-z_v) J_0(u(z-z_v)) dz_v.
\end{equation}
Using the integral identity of the Bessel function $ \int_{0}^{x} v J_0(v)\, dv=x J_1(x)$, we can write equation (\ref{bessel1}) as follows
\begin{equation}\label{bessel2}
\int_{x}^{z_{max}}G_{\lambda,\mu}(z_v)dz_v=\frac{1}{u}\, \int_{x}^{z_{max}}\widehat{f}_{\lambda,\mu}(z)(z-x) J_1(u(z-x)) dz.
\end{equation}

\noindent Adding equations (\ref{bessel2}) and (\ref{diff}) we get
\begin{equation}\label{final1}
\displaystyle {{\mathcal{{H}}}\left(\int_{x}^{z_{max}}G_{\lambda,\mu}(z_v)dz_v\right)=\,\int_{x}^{z_{max}}\widehat{f}_{\lambda,\mu}(z) J_0\left(u(z-x)\right)\,  dz}.
\end{equation}
\noindent We again apply the operator ${\mathcal{{H}}}$ to equation (\ref{final1}) to obtain

\begin{eqnarray*}
{\mathcal{{H}}}^2 \left(\int_{x}^{z_{max}}G_{\lambda,\mu}(z_v)dz_v\right)=-\widehat{f}'_{\lambda,\mu}(x)+u^2 \,\int_{x}^{z_{max}}\widehat{f}_{\lambda,\mu}(z) [J_0\left(u(z-x)\right)\\
\hspace{+5.5cm} + J_0''\left(u(z-x)\right)\, ] dz.
\end{eqnarray*}

\noindent Then we integrate our equation as follows
\begin{eqnarray}\label{long}
\hspace{-1.5cm}\int_{t}^{z_{max}} J_0\left(u(t-x)\right){\mathcal{{H}}}^2 \left(\int_{x}^{z_{max}}G_{\lambda,\mu}(z_v)dz_v\right)\,dx=-\int_{t}^{z_{max}} J_0\left(u(t-x)\right)\widehat{f}'_{\lambda,\mu}(x)dx\nonumber \\
 \hspace{-1.5cm}+ u^2 \int_{t}^{z_{max}} J_0\left(u(t-x)\right)\int_{x}^{z_{max}}\widehat{f}_{\lambda,\mu}(z)[ J_0\left(u(z-x)\right)+ J_0''\left(u(z-x)\right)]dz\,dx.
\end{eqnarray}
\noindent  Let us simplify the RHS of equation (\ref{long}). First, we simplify the first term of the RHS by integration by parts
\[\int_{t}^{z_{max}} J_0\left(u(t-x)\right)\widehat{f}'_{\lambda,\mu}(x)dx=\widehat{f}_{\lambda,\mu}(t)-u\int_{t}^{z_{max}} J_1\left(u(t-x)\right)\widehat{f}_{\lambda,\mu}(x)dx.
\]
\noindent We used the conditions imposed upon the function $\widehat{f}_{\lambda,\mu}(z_{max})=0$. Then we simplify the second term of the RHS by interchanging the order of integration
\begin{eqnarray*}
 \int_{t}^{z_{max}} J_0\left(u(t-x)\right)\int_{x}^{z_{max}}\widehat{f}_{\lambda,\mu}(z)[ J_0\left(u(z-x)\right)+ J_0''\left(u(z-x)\right)]dz\,dx \\
 =\int_{t}^{z_{max}}\widehat{f}_{\lambda,\mu}(z)dz\,\int_{t}^{z}J_0\left(u(t-x)\right)[ J_0\left(u(z-x)\right)+ J_0''\left(u(z-x)\right)]dx.
\end{eqnarray*}
\noindent Then we use the recurrence formula of Bessel function $J_0''(x)+J_0(x)=-J_1(x)/x$ and the formula
$\displaystyle {\int_{0}^{z}J_1(x) J_0(z-x) \frac{dx}{x}=J_1(z)}$ (for more details about Bessel functions see \cite{bessel} page 671). After simplification, we obtain that for $(\lambda,\mu)\ne(0,0)$

 \[ \widehat{f}_{\lambda,\mu}(t)=\displaystyle{\int_{t}^{z_{max}} J_0\left(u(t-x)\right){\mathcal{{H}}}^2 \left(\int_{x}^{z_{max}}G_{\lambda,\mu}(z_v)dz_v\right)\,dx}.\]
Since $f$ is compactly supported, its Fourier transform at $(\lambda,\mu)=(0,0)$ is uniquely defined by continuity.

\hfill $\Box$

\section{Numerical implementation}
To validate our inversion formulae in 2D we run a few simple numerical simulations. Here we present the results of these experiments for smooth phantoms supported within a square $\Omega=[-1,1]\times [-1,1]$ and defined by functions of the form
\begin{eqnarray*}
\hspace{-2cm}f(x,y)=\left\{
\begin{array}{ll}
exp\left\{\displaystyle\frac{-r^2}{r^2-\left[(x-x_c)^2+(y-y_c)^2\right]}\right\} & \mbox{ if  }\, (x-x_c)^2+(y-y_c)^2  < r^2 , \\ 0  & \mbox{ otherwise . } \end{array} \right.
\end{eqnarray*}

In all experiments we used a fixed opening angle $\beta=\pi/8$. The V-line Radon transform was computed numerically from an $N\times N$
discretized version of the function $f$. The differentiation operator in the inversion formula (\ref{result}) was approximated by the standard first order forward difference formula
\[ \frac{\partial}{\partial z} g(x,z)\approx \frac{ g(x,z+\triangle z)- g(x,z)}{\triangle z},\]
and the integration was done using the trapezoidal rule.

\begin{figure}[h]
\begin{center}
{\epsfig{figure=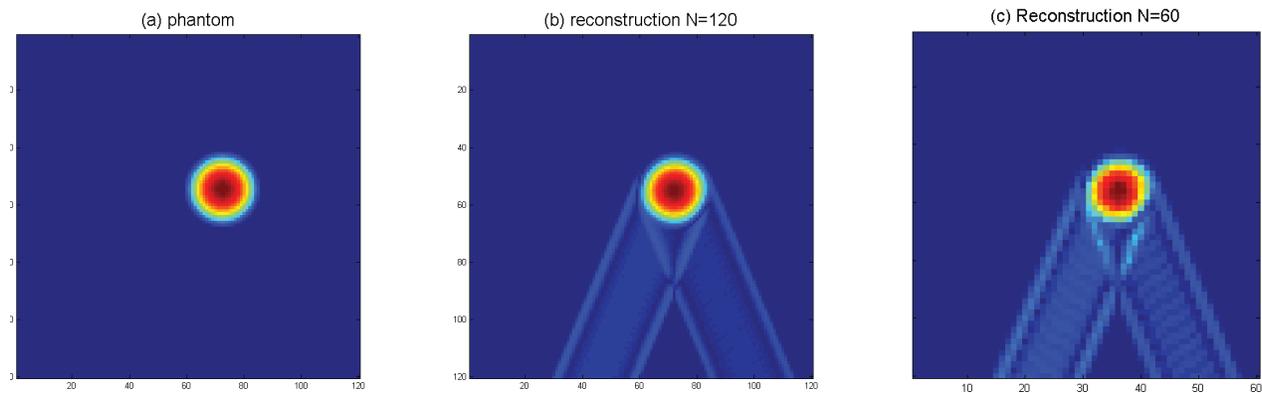,height=5in}}
\end{center}
\vspace{-4.5cm}
\caption{Numerical results: (a) phantom, (b) reconstruction using $N=120$ and (c) reconstruction using $N=60$}  \label{numerical1}
\end{figure}

Figure \ref{numerical1} shows the results for a phantom with $r=0.25$ and the center $(x_c,y_c)=(0.2,0.1)$ using two different values of discretization $N=60$ and $N=120$ pixels.

\begin{figure}[h]
\begin{center}
{\epsfig{figure=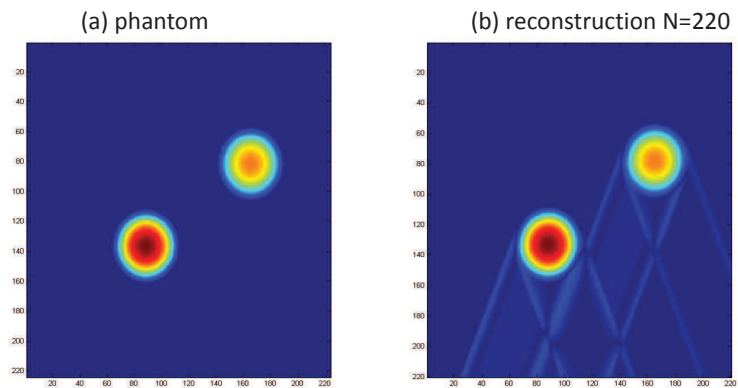,height=4in}}
\end{center}
\vspace{-2.5cm}
\caption{Image reconstruction for N=220 using $\beta=\pi/8$} \label{circles}
\end{figure}

We also used a phantom consisting of two circles of
varying centers and intensities as shown in table
\ref{table}. The images were constructed with 220$\times$ 220 pixels (figure \ref{circles}).
\begin{table}[htbp]

  \centering

  \caption{Parameters for the phantom image}\label{table}

\begin{tabular}{|c|c|c|c|}

\hline

Circle&Coordinates of the center$(x_c,y_c)$&intensity\\

\hline

1&(0.5,0.3)&3\\

\hline

2&(-0.2,-0.2)&4 \\

\hline

\end{tabular}

  \label{tableintro:1}

\end{table}

\section{Additional remarks}
\begin{enumerate}
\item The smoothness and decay conditions for $f$ in both inversion formulae (for 2D and 3D) are not optimized. The formulae may hold with weaker requirements, e.g. for $f$ in Schwartz space, or compactly supported functions that have only a few orders of smooth derivatives.
\item The main results of this paper have potential to be generalized to dimensions higher than $n=3$. The authors plan to address this issue in future work.
\item The numerical implementation of the 3D inversion formula presented here is a challenging task in itself and the authors plan to address this in future work.
\end{enumerate}

\ack
The first author was supported by the American University of Sharjah (AUS) research grant FRG3. The work of the second author was supported in part by US NSF Grant DMS 1109417.

A part of this paper was written during the second author's one semester long visit to the American University of Armenia (AUA). He thanks the administration and staff of AUA for their support and stimulating environment.

\end{document}